\documentclass{amsart}
\usepackage{amssymb}
\usepackage{amscd}
\usepackage{verbatim}
\usepackage{epsfig}
\begin{document}

\newcommand\RR{\mathbb{R}}
\newcommand\CC{\mathbb{C}}
\newcommand\Lap{\Delta}
\newcommand\pa{\partial}
\newcommand\ep{\epsilon}
\renewcommand\a{\alpha}
\renewcommand\b{\beta}
\renewcommand\l{\lambda}
\newcommand\symbint{\sigma_{\operatorname{int}}}
\newcommand\symbb{\sigma_{\pa}}
\newcommand\Dprime{{\mathcal D}'}
\newcommand\td{\tilde}
\def\ang#1#2{\big\langle #1 , #2 \big\rangle}
\newcommand\ol{\overline}
\newcommand\hz{\hat z}
\newcommand\Nu{{\mathcal V}}
\renewcommand\sp{\operatorname{sp}}
\newcommand\lb{\operatorname{lb}}
\newcommand\rb{\operatorname{rb}}
\newcommand\bfc{\operatorname{bf}}
\newcommand\ul{\underline}
\newcommand\Ls{L^{\sharp}}
\newcommand\Lt{\tilde L}
\newcommand\tr{\operatorname{tr}}
\newcommand\rtr{\operatorname{r-tr}}
\newcommand\spl{s(\l)}
\newcommand\Lapext{\Delta_{\Omega}}
\newcommand\Lapplane{\Delta_{\RR^2}}
\newcommand\Lapint{\Delta_{\mathcal{O}}}
\newcommand\Lapsum{\Delta_{\oplus}}
\newcommand\Lapplus{\Delta_{\oplus}}
\newcommand\Lapprime{\Delta_{\Omega'}}
\newcommand\CI{C^\infty}
\newcommand\Vol{\operatorname{Vol}}
\newtheorem{lemma}{Lemma}[section]
\newtheorem{prop}[lemma]{Proposition}
\newtheorem{thm}[lemma]{Theorem}
\newtheorem{cor}[lemma]{Corollary}
\newtheorem{result}[lemma]{Result}
\newtheorem*{thm*}{Theorem}
\numberwithin{equation}{section}
\numberwithin{figure}{section}
\theoremstyle{remark}
\newtheorem{rem}[lemma]{Remark}
\theoremstyle{definition}
\newtheorem{defn}[lemma]{Definition}
\newcommand\downto{\to}
\newcommand\upto{\to}
\newcommand\Pdir{P_{\operatorname{dir}}}
\newcommand\Ptr{P_{\operatorname{tr}}}
\newcommand\Tdir{T_{\operatorname{dir}}}
\newcommand\Ttr{T_{\operatorname{tr}}}
\newcommand\Rplus{{\RR_+}}
\newcommand\ttr{\tilde {\operatorname{tr}}}
\newcommand\rone{r_1}
\newcommand\rtwo{r_2}
\newcommand\rthr{r_3}
\newcommand\HR{\operatorname{HR-}\!\!\!}
\newcommand\ilg{\operatorname{ilg}}
\newcommand\zbar{\overline{z}}
\newcommand\intcirc{\int_{S^1}}
\newcommand\intdisc{\int_{D}}
\newcommand\phit{\phi_t}
\newcommand\op{\operatorname{op}}
\newcommand\just{\mskip 1000mu minus998mu}
\newcommand\obst{\mathcal{O}}
\newcommand\bdy{\partial \mathcal{O}}
\newcommand\half{\frac{1}{2}}
\newcommand\uol{u_{\omega, \lambda}}
\newcommand\Id{\operatorname{Id}}
\newcommand\Area{\operatorname{Area}}
\newcommand\loc{\operatorname{loc}}

\title[Determinants of Laplacians in exterior domains]{Determinants of 
Laplacians in exterior domains}

\author{Andrew Hassell}
\address{Centre for Mathematics and its Applications, Australian National
University, Canberra ACT 0200, Australia}
\email{hassell@maths.anu.edu.au}
\author{Steve Zelditch}
\address{Department of Mathematics, Johns Hopkins University, 
Baltimore, MD 21218, USA}
\email{zel@math.jhu.edu }

\thanks{Research of the first author supported by Australian Research Council;
research of the second author 
partially supported by  NSF grant \#DMS-9500491}

\begin{abstract} We consider classes of simply connected planar domains which
are {\it isophasal}, ie, have the same scattering phase $s(\l)$ for
all
$\l > 0$. This is a
scattering-theoretic analogue of isospectral domains. Using the heat invariants 
and the determinant of the Laplacian, Osgood, Phillips and
Sarnak showed that each isospectral class is sequentially compact in a natural 
$\CI$ topology.
This followed earlier work of Melrose who showed that the set of 
curvature functions $k(s)$ is compact in $\CI$. 

In this paper, we show sequential compactness of each isophasal class
of domains. To do this
we define the determinant of the exterior Laplacian and use it together with 
the heat invariants (the heat invariants and the determinant being 
isophasal invariants). We show that the
determinant of the interior and exterior Laplacians satisfy a 
Burghelea-Friedlander-Kappeler type surgery formula. This allows a reduction to
a problem on bounded domains for which the methods of Osgood, Phillips
and Sarnak can be adapted. 
\end{abstract}
\maketitle


\section{Introduction}

\subsection{The isospectral problem}

\begin{figure}\label{pinch}\centering
\epsfig{file=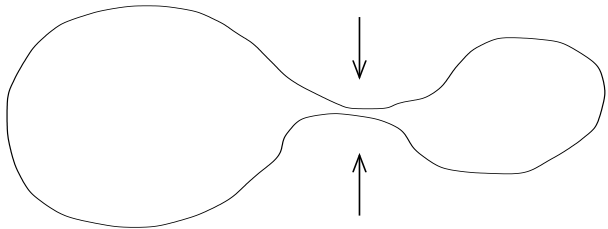,width=8cm,height=3cm}
\caption{Pinching off without blowup of curvature}
\end{figure}

In this paper, we consider a scattering-theoretic version of the famous
question `Can one hear the shape of a drum?' posed by M. Kac \cite{Kac}. In
mathematical terms the question is whether a planar domain $\mathcal{O}$
is determined up to isometry by its Laplace spectrum (with Dirichlet
boundary conditions, say), where the spectrum $0 < \l_1^2 \leq \l^2_2 
\dots$ is counted with multiplicity. The answer to this question is
known to be negative \cite{GWW} (though there are some positive results
for restricted classes of domains \cite{Zel}). 
In view of this, it is reasonable to ask how
`small' is the set of domains isospectral to a given domain. One
way to make this precise is to ask whether the isospectral class is
compact in some topology on domains. Melrose \cite{Mel1} showed that this is
the case, where the topology is taken as the $\CI$ topology on the
curvature function $k(s) : s \in [0, L]$ of the boundary of the domain
($L$ is fixed over the isospectral class, as discussed below). A
disadvantage of this topology is that it does not exclude the
possibility of a sequence of isospectral domains pinching off (see 
figure~\ref{pinch}).

This result was proved using the `heat invariants'. The heat invariants
of a domain $\mathcal{O}$ are coefficients in an expansion of the heat
trace $e_{\mathcal{O}}(t)$ as $t \to 0$. Since 
$$
e_{\mathcal{O}}(t) \equiv \tr e^{-t\Lapint} = \sum_{j=1}^\infty
e^{-t\l_j^2},
$$
the heat trace is a spectral invariant. It has a well known asymptotic
expansion
\begin{equation}
e_{\mathcal{O}}(t) \sim \sum_{j=-2}^\infty a_j t^{-j/2}, \quad t \to 0
\label{ae}
\end{equation}
as $t$ tends to zero, with the heat invariants $a_j$ `local', that is,
integrals over the domain $\mathcal{O}$ or the boundary $H = 
\pa \mathcal{O}$ of locally defined
geometric quantities. The first few are
\begin{equation}\begin{aligned}
a_{-2} &= \frac{\text{area}(\mathcal{O})}{4\pi} = \frac1{4\pi}\int_{\mathcal{O}}
1 \\
a_{-1} &= c_{-1} \text{ length}(\bdy) = c_{-1} \int_{\bdy} 1 \, ds \\
a_{0} &= c_0 \chi(\obst) = c_{0} \int_{\bdy} k(s) \, ds \\
a_1 &= c_1 \int_{\bdy} k^2(s) \, ds
\end{aligned}\end{equation}
with $c_i \neq 0$. Melrose showed that
\begin{equation}
a_{2l-1} = c_l \int_{\bdy} (k^{(l)}(s))^2 + p_l(k(s), \dots k^{(l-1)}(s))
\, ds \quad c_l \neq 0
\label{odd-heat-inv}
\end{equation}
where $p_l$ is a polynomial; this is the main step in his
result. 

Notice that the first two heat invariants give $A = $ area$(\obst)$ and
$L = $ length$(\bdy)$ as spectral invariants. Consideration of the
isoperimetric quotient shows that the disc of radius $r$, $D_r$, is determined
by its spectrum, an observation that perhaps led to Kac's question. 

Osgood, Phillips and Sarnak (abbreviated OPS from here on) considered the
question of $\CI$ compactness from a different point of view
\cite{OPS1}, \cite{OPS2}. They regarded
a planar domain as the image of the unit disc $D$ under a conformal map $F$.
The metric on the domain is then isometric to $e^{2\phi} g_0$ where
$g_0$ is the flat metric on the disc, and $\phi = \log |F'|$ is a harmonic
function. Thus $\phi$ is determined by its boundary values $\phi \restriction
\pa D$. Given a harmonic function $\phi$, one can find the
corresponding domain $F(D)$, which is a flat planar domain 
(possibly self-overlapping). 
OPS showed that isospectral classes are compact in the $\CI$
topology of $\phi$ restricted to $\pa D$. This result excludes degenerations
of the form illustrated in figure~\ref{pinch} since derivatives of $F$ must
blow up under such a degeneration. (Melrose also has an argument ruling out
such degenerations using the first positive singularity of the wave trace
\cite{Mel2}.)

Osgood, Phillips and Sarnak used the heat invariants, plus one other
invariant, the determinant of the Laplacian (described below), 
to deduce their result. They exploit a formula, due to
Polyakov \cite{Poly} and Alvarez \cite{A}, 
expressing the determinant in terms of the function $\phi$:
\begin{equation}
\log \det \Lapint = \frac{1}{12\pi} \intcirc \phi \pa_n \phi 
- \frac{1}{6\pi} \intcirc \phi + \log \det \Delta_{D}.
\end{equation}
The formula is remarkable since the first term on the right is nonnegative
and almost the square of the Sobolev $\half$ norm of $\phi$. 

\subsection{An analogous problem for the exterior domain}\label{exterior} 
In this paper we
are interested in the exterior Laplacian. Let $\Omega = \RR^2 \setminus
\obst$, the exterior of an obstacle, and let $\Lapext$ be the Laplacian
on $L^2(\Omega)$ with domain $H^2(\Omega) \cap H^1_0(\Omega)$ (ie, the
Dirichlet Laplacian). The operator $\Lapext$ is self-adjoint with continuous
spectrum on $[0, \infty)$ (\cite{Taylor}, chapter 8). 
We wish to formulate a problem about exterior
domains that is analogous to the isospectral problem. 

To do this, we observe that the isospectral condition may be expressed in
terms of the counting function
$$
N_\obst (\l) = \text{ number of eigenvalues of } \Lapint \leq \l^2.
$$
To say that two domains have the same spectrum, counted with multiplicity,
is equivalent to saying that they have the same counting function $N(\l)$,
and the problem considered by OPS is to show that the class of domains with
a fixed $N(\l)$ are compact in some natural topology. 

For exterior domains, $\Omega = \RR^2 \setminus \obst$, it is known that
$1/2\pi$ times
the scattering phase $\spl$ is analogous to the counting function. The
usual definition of the scattering phase is
$$
\spl = -i \log \det S_\Omega (\l),
$$
where $S_\Omega (\l)$ is the scattering matrix (see \cite{gst}). For our
purposes, it is more illuminating to note that the difference between
the spectral projection $E_\Omega (\l)$ on the interval $(-\infty, \l)$
for $\Lapext$, and the corresponding spectral projection $E_0(\l)$ for
$\Lapplane$, is trace class in a distributional sense, with
\begin{equation}
\tr \int_0^\infty \phi'(\sigma) \big( E_\Omega (\sigma) - E_0(\sigma) \big) 
d\sigma =
\tr \big( \phi(\Lapext) \oplus 0 - \phi(\Lapplane) \big) = 
\int_0^\infty \phi'(\sigma) \frac{s(\sqrt{\sigma})}{2\pi} \, d\sigma 
\label{funct-calc}\end{equation}
for any $\phi \in C_c^\infty(\RR)$ (see \cite{JK}; the normalization
of their scattering phase $\theta(\l)$ is minus one-half of our
$s(\l)$ --- cf remark 1 of their introduction). 
Thus $-s(\l)/2\pi$ is a regularized trace of the spectral measure.
Since
$$
N_\obst (\sqrt{\sigma}) = \tr E_\obst (\sigma),
$$
the analogy between the counting function and the scattering phase is
clear. Let us say that two obstacles are {\it isophasal} if they have
the same scattering phase.

A strong indication that it might be possible to use the scattering phase
to prove compactness results about isophasal classes of domains comes from
noting that the formula \eqref{funct-calc} holds also for 
$\phi(\sigma) = e^{-\sigma t}$. Thus the regularized trace of the heat kernel is
given in terms of the scattering phase by
\begin{equation}
\rtr e^{-t\Lapext} \equiv 
\tr \big( e^{-t\Lapext} - e^{-t\Lapplane} \big) = -\frac{t}{\pi} \int_0^\infty
\spl e^{-\l^2 t} \l d\l .
\label{r-tr-sc-phase}
\end{equation}
However, direct construction of a parametrix for the heat kernel of 
$\Lapext$ near $t=0$ shows that the regularized trace has an asymptotic
expansion of the form \eqref{ae} with the {\it same } coefficients
(up to changes of sign). Thus we immediately get Melrose's result for
$\CI$ compactness of the curvature function of the boundary. The question
is then whether this can be improved, OPS-style, to a result of $\CI$
compactness of the domain. It is very natural to look for an analogue of
the determinant of the exterior operator in order to do this.

\subsection{Determinants and surgery formulae}\label{dets} We begin by recalling the
definition of the determinant. Let $A$ be a strictly positive elliptic
$m$th order differential
operator on a bounded domain of dimension $n$
(compact manifold, possibly with boundary). 
Then $A$ has positive, discrete spectrum $0 < \mu_1 \leq \mu_2 \dots \to
\infty $. The determinant of $A$ is defined in terms of the zeta function,
$\zeta(s)$. The zeta function is defined by
\begin{equation}
\zeta(s) = \sum_{j=1}^\infty \mu_j^{-s} = \frac1{\Gamma(s)} \int_0^\infty
t^s \tr e^{-tA} \frac{dt}{t}
\end{equation}
in the region of absolute convergence, $\Re s > n/m$. Since the heat 
trace has an expansion \cite{Seeley}
\begin{equation}
\tr e^{-tA} \sim \sum_{j = -n}^\infty t^{-j/m} b_j, \quad t \to 0,
\label{A-heat-exp}\end{equation}
it follows that $\int t^{s-1} \tr e^{-tA} dt$ continues meromorphically to
the complex plane with at most simple poles at $s = -j/m, j \geq -n$. The
factor $\Gamma(s)^{-1}$ vanishes at $s=0$ ensuring that the zeta function
is regular at $s=0$. The determinant of $A$ is then defined by
$$
\log \det A = -\zeta'(0).
$$
If $A$ is not strictly positive, that is, has a zero eigenvalue, then the
determinant is defined to be zero. However, it is usually of interest to
look instead at the modified determinant, $\det {}' A$. This is defined
by defining the zeta function using only the nonzero eigenvalues of $A$, and
then taking $\log \det {}' A = -\zeta'(0)$. An equivalent definition is that
the modified determinant of $A$ is the determinant of $A + \Pi_0$, where
$\Pi_0$ is orthogonal projection onto the null space of $A$. 

If $A$ is pseudodifferential, then the definition above does not make sense in
general, since 
the heat trace of $A$ may have log terms (terms of the form $t^{j/m} \log t$)
as $t \to 0$, and then the zeta function may have a pole at $s=0$. However,
in the case of interest in this paper --- the Neumann jump operator
(see Definition~\ref{Neumann-jump}) --- one can rule this out and then
the log determinant is defined just as for differential operators. 

The terminology `determinant' is justified by the fact that if $A$ were an
operator on a finite dimensional space, and therefore had
a finite number of eigenvalues, then we would have
$$
-\zeta'(0) = - \sum (-\log \mu_j) \mu_j^{-s} |_{s=0} = \sum \log \mu_j = \log 
\det A.
$$

\

The log determinant is a non-local quantity; that is, it cannot be
written as the integral over $\obst$ or $\bdy$ of locally-defined
geometric quantities \cite{Ray}. However, it behaves in many situations as
a `quasi-local' quantity, in the following sense: when a localized
perturbation is made in the operator, one can often find a formula for the
change of the log determinant which involves only the perturbation. 
An example, which is highly relevant to this paper, is the Mayer-Vietoris
type surgery formula for the log determinant proved by Burghelea, Friedlander
and Kappeler \cite{BFK} (henceforth BFK). 
Suppose that $A$ is an elliptic partial
differential operator on a compact manifold $M$, and $H$ is a hypersurface,
with $\tilde M$ the manifold with boundary obtained by cutting $M$ at $H$.
Let $B$ be an elliptic boundary condition for $A$ on the boundary of 
$\tilde M$. BFK 
found a formula for $\log \det A - \log \det (A, B)$ in terms
of the log det of a 
pseudodifferential operator $R$ on $H$ and other data defined on $H$. 
In the particular case of the Laplacian on a two-dimensional manifold $M$,
with $H$ a curve dividing $M$ into two components $M_1$ and $M_2$,
and $B$ the Dirichlet boundary condition, they showed that
\begin{equation}
\log \det \Delta_M - \log \det (\Delta_{M_1}, B) - \log 
\det (\Delta_{M_2}, B) = \log \det R + \log a - \log l,
\label{BFK}
\end{equation}
where $R$ is the Neumann jump operator (see Definition~\ref{Neumann-jump}), 
$a$ is the area of $M$ and $l$
is the length of $H$. 

In this paper, we look at the exterior and interior Dirichlet Laplacians
for an obstacle $\obst$ from this point of view. Thus, we consider the boundary
$H$ of $\obst$ to be a cutting of the manifold $\RR^2$ and look for a
BFK-type surgery formula. Of course, one problem is that $\RR^2$ is
unbounded, so $\Lapext$ has continuous spectrum and its log determinant
is not defined. This is the topic of the next section. The main theorem
is Theorem~\ref{surg-formula}, which gives a surgery formula for this
regularized log determinant very similar to \eqref{BFK}. The proof of
this theorem is the subject of the third section. 

\subsection{Compactness of isophasal sets}\label{isophasal} 
In the fourth section we show
that each class of isophasal sets is compact in a natural $\CI$ topology
(Theorem~\ref{compactness}).
First we must specify the topology on domains. Following Osgood,
Phillips and Sarnak, we define a sequential topology, ie, we specify
the convergent sequences rather than the open sets. This is
appropriate since our goal is to prove sequential compactness.

It is inconvenient to deal
with unbounded domains, so we pass to the inversion $\Omega^I$ of
$\Omega$. We will say that a sequence of exterior domains $\Omega_i$
converges in the $\CI$ topology if there are Euclidean motions $E_k$
of the plane such that 

(i) the closure of $E_k \Omega_k$ does not contain the origin;

(ii) the sequence $\obst _k$ of inversions of $E_k \Omega_k$ about
the origin converges in the sense that there are conformal maps $F_k$
from the disc to $\obst _k$, with $|F_k'|$ never zero, which converge 
in the $H^s$ topology for all $s$. 

The notion of convergence of $\obst _k$ in condition (ii) is stronger
than convergence in the OPS topology, which would say that there are
Euclidean motions $\td E_k$ and conformal maps $F_k$ from the disc to
$\td E_k \obst _k$ that converge in $H^s$ for all $s$. It is important
that the group of Euclidean motions is allowed to act on the domains
$\Omega_k$ and not on the $\obst _k$.

To prove Theorem~\ref{compactness}, 
we use the conformal equivariance of the Laplacian
in two dimensions to compactify the problem. That is, we consider a metric
$g$ on $\RR^2$ which is a conformal multiple of the flat metric, so that
infinity is compactified to a point. Comparing our surgery formula 
to that of BFK on the compactified space, we show that the Laplacian
on a certain bounded domain of $S^2$ (depending on the obstacle)
has fixed determinant, as the obstacle
ranges over an isophasal set. This allows us to adapt the argument of
OPS to obtain the result. 

\section{Determinant of the exterior operator}
In this section we shall define the modified determinant of the exterior
Laplacian, and state the main theorem. 

The determinant is usually defined in terms of the zeta function, which in turn is
defined in terms of the trace of the heat kernel. In the case of the exterior
Laplacian, the heat kernel is certainly not trace class, since it has continuous
spectrum. However, as discussed in section~\ref{exterior}, the difference between
the exterior heat operator and the free heat operator is trace class for every
$t$ (see \cite{BK} and \cite{JK}); 
we will denote this trace by $\rtr e^{-t\Lapext}$, and call it the
regularized heat trace. Thus, the obvious candidate for the zeta function is 
\begin{equation}
\zeta_\Omega (s) \ `= \text{'} \ 
\frac1{\Gamma(s)} \int_0^\infty t^s \rtr e^{-t \Lapext} \frac{dt}{t}.
\label{zeta-heattrace}
\end{equation}

Unfortunately, this integral does not converge for any value of $s$. 
To deal with
this we break up the regularized heat trace into two pieces. Let $\chi$ be a smooth
function that it identically one near $\l = 0$ and identically zero for $\l > 1$.
Then, by \eqref{r-tr-sc-phase}, $\rtr e^{-t\Lapext}
= e_1(t) + e_2(t)$, where
\begin{equation}
e_1(t) = -\frac{t}{\pi} \int_0^\infty \chi(\l) \spl e^{-\l^2 t} \l d\l 
\label{e1}
\end{equation}
and
\begin{equation}
e_2(t) = -\frac{t}{\pi} \int_0^\infty (1 - \chi(\l)) \spl e^{-\l^2 t} \l d\l .
\label{e2}
\end{equation}
We write $\zeta_\Omega (s) = \zeta_{\Omega, 1}(s) + \zeta_{\Omega_2}(s)$ for the
corresponding decomposition of the zeta function. 
Then $e_2(t)$ is exponentially decreasing at infinity. On the other hand, since
$\rtr$ has the usual asymptotic expansion at $t=0$, and $e_1(t)$ is smooth at
$t=0$, we see that $e_2(t)$ has an expansion of the form \eqref{A-heat-exp}.
Thus, 
\begin{equation}
\zeta_{\Omega,2} =  \frac1{\Gamma(s)} \int_0^\infty t^s e_2(t) \frac{dt}{t}
\label{zeta-2-heat}
\end{equation}
continues meromorphically to the entire plane with no
pole at $s=0$ by the usual argument. The other part, $\zeta_{\Omega, 1}(s)$ may
be directly expressed in terms of the scattering phase by the formula
\begin{equation}
\zeta_{\Omega,1} (s) = -\frac{s}{\pi} \int_0^\infty \l^{-2s} \spl
\chi(\l) \frac{d\l}{\l}.
\label{zeta-1-sp}
\end{equation}
To understand this integral we need the asymptotics of $\spl$ as $\l \to 0$. 
Following \cite{HMM}, 
let us define, for this and other purposes, the function ilg.

\begin{defn}\label{ilg} The function $\ilg \l$ is defined to be
$$
\ilg \l = \frac{1}{\log (1/\l)};
$$
it goes to zero as $\l \to 0$, but slower than any positive power of $\l$.
\end{defn}

\begin{lemma}\label{sc-phase-zero} 
For any obstacle $\mathcal{O}$, the scattering phase satisfies
\begin{equation}
\spl = \pi \ilg \l + O((\ilg \l)^2), \quad \l \to 0.
\label{ltozero}
\end{equation}
Here, the $O((\ilg \l)^2)$ term is
uniform over each isophasal class.
\end{lemma}

\begin{proof} In the appendix we compute $\spl$ for a disc of radius $1$; the result is
$$
s_{D_1}(\l) = \pi \ilg \l + O((\ilg \l)^2), \quad \l \to 0. 
$$
The scattering phase for a disc of radius $r$ is then $s_{D_r}(\l) = s_{D_1}(r^2 \l)$.
Thus, for a disc of radius $r$, we also have 
$$
| s_{D_r}(\l) - \pi \ilg \l | \leq C(r) (\ilg \l)^2.
$$
Using the first two heat invariants, the area and perimeter are constant on an
isophasal class of domains; hence, by Lemma~\ref{inradius} 
the inradius and circumradius are
uniformly bounded below and above. Thus, we can sandwich any domain in an isophasal
class between fixed discs $D_r$ and $D_R$. 
By \cite{HR}, $s(\l)$ is monotonic in the domain, so we obtain
\eqref{ltozero}. 
\end{proof}

Substituting this expansion into \eqref{zeta-1-sp}, we see that the zeta function is
meromorphic in the half plane $\Re s < 0$, but {\it not} in any neighbourhood of
$s=0$. To see this, consider the function
$$
g(s) = \int_0^\infty \l^{-2s} \ilg \l \chi(\l) \frac{d\l}{\l}.
$$
By differentiating once in $s$, it is not hard to show that $g(s)$ is equal
to $-\log (-s)$ plus a smooth function as $s \uparrow 0$. Thus, 
the zeta function has an expansion of the form
\begin{equation}
\zeta_\Omega (s) = a_0 - s \log (-s) + a_2 s + O(s^2 \log s), \quad s 
\uparrow 0.
\label{zeta-exp}
\end{equation}
This allows us to make 

\begin{defn}\label{detext} The logarithm of the
determinant of the exterior Laplacian is defined to be
$$
\log \det {}' \Lapext = -a_2,
$$
where $a_2$ is the coefficient of $s$ in the expansion \eqref{zeta-exp}.
\end{defn}

For future use, we observe here that if we consider the operator $\Lapext + \mu$
instead of $\Lapext$, with $\mu > 0$, then the zeta function is given by
\begin{equation}
\zeta_{\Omega, \mu} (s) = \int_0^\infty t^s \rtr e^{-t\Lapext} e^{-\mu t} \frac{dt}{t};
\label{zeta-mu-heat}
\end{equation}
the integral is now defined for $\Re s > 1$ and continues meromorphically to the
complex plane with no pole at $s=0$, since the exponential factor $e^{-\mu t}$ makes
the integral convergent at infinity for any $s$. However, it is useful to write the
zeta function in the same way as for the case $\mu = 0$:
\begin{equation}
\zeta_{\Omega, \mu} (s) = 
-\frac{s}{\pi} \int_0^\infty (\l^2 + \mu)^{-s-1} \chi(\l)
\spl \l d\l + \frac1{\Gamma(s)} \int_0^\infty t^s e_2(t) e^{-\mu t} \frac{dt}{t}.
\label{zeta-mu-chi}
\end{equation}
The determinant is then defined in the usual way. 

\begin{defn} For $\mu > 0$, the logarithm of the
determinant of $\Lapext + \mu$ is defined by
$$
\log \det (\Lapext + \mu) = - \zeta_{\Omega, \mu}'(0).
$$
\end{defn}

{\it Remark. } We write $\det {}'$ instead of $\det$ in definition~\ref{detext}
because it is more similar to the modified determinant described in 
section~\ref{dets} than the determinant --- this becomes clear in the calculation
of section~\ref{mutozero}. 

\vskip 10pt

It is not yet clear that the quantity in Definition~\ref{detext} merits the
term `determinant'. 
We believe
that the following theorem justifies the definition --- 
compare with equation \eqref{BFK}. 
First we give a formal
definition of the Neumann jump operator $R$. Before stating it, we observe
that for any $\mu \geq 0$, and given any continuous function $f$ on $H$,
there is a unique bounded extension $u$ of $f$ to $\Omega$ satisfying
$(\Lapext + \mu) u = 0$. 

\begin{defn}\label{Neumann-jump} 
The Neumann jump operator $R$ for the obstacle $\obst$ is
the operator 
$$
f \mapsto \pa_\nu u_1 - \pa_\nu u_2,
$$
where $f \in C^0(H)$, $u_1$, respectively $u_2$ are the bounded extensions of
$f$ to $\obst$, respectively $\Omega$ satisfying $\Delta u_i = 0$,
and $\nu$ is the outward normal. (This choice of normal means that
$R$ is a nonnegative operator.)
The operator $R(\mu)$ is
defined similarly, replacing $\Delta$ with $\Delta + \mu$. The operators
$R$ and $R(\mu)$ are pseudodifferential operators of order $1$, and
for $\mu > 0$, $R(\mu)$ is strictly positive. 
\end{defn}

\begin{thm}\label{surg-formula} 
The following formula holds:
\begin{equation}
\log \det {}' \Lapext + \log \det \Lapint + \log \det {}' R = 
\gamma + \log \frac{L}{\pi}
\end{equation}
where $\gamma$ is Euler's constant. 
\end{thm}

The proof of this theorem is the subject of the next section.

\section{Proof of the surgery formula}
In this section we prove Theorem~\ref{surg-formula}. We follow closely the scheme of
BFK's proof. Thus, the proof consists of three steps. The first step is to
establish the variational formula
\begin{equation}
\frac{d}{d\mu} \big( \log \det  (\Lapext + \mu) + \log \det (\Lapint + \mu) +
\log \det R(\mu) \big) = 0,
\label{var-formula}
\end{equation}
for $\mu > 0$,
where $R(\mu)$ is the Neumann jump operator. This calculation 
was first done by Forman \cite{Forman} and the
proof given here is almost identical, but it is written out in full for the reader's
convenience.
Thus, integrating \eqref{var-formula}, we find that
\begin{equation}
\log \det  (\Lapext + \mu) + \log \det (\Lapint + \mu) + \log \det R(\mu) = C.
\label{C}
\end{equation}

In the second step, we show that $C = 0$. To do this, 
we send $\mu$ to infinity. Then each of the
log determinants has an asymptotic expansion in $\mu$, with local coefficients.
Clearly the coefficients of each term must agree on the left and right hand
side of \eqref{C} so if we know the coefficient of the constant term for each
log determinant, then we deduce the value of $C$.
It turns out that the constant term in
the expansion for each log determinant is zero, so $C = 0$.

The third step is to consider the limit $\mu \to 0$. Here 
we prove the following asymptotic expansions:
\begin{align}
\log \det (\Lapext + \mu) &= \log \log \mu^{-1/2} + \log \det {}' \Lapext
- \gamma - \log 2 + o(1), \ \mu \to 0; \label{ext-asympt-0}\\
\log \det (\Lapint + \mu) &= \log \det \Lapint + o(1), \ \mu \to 0; 
\label{int-asympt-0} \\
\log \det R(\mu) &= -\log \log \mu^{-1/2} +
\log \det {}' R - \log \frac{L}{2\pi} + o(1), \ \mu \to 0. \label{R-asympt-0}
\end{align}
It is easy to see that  Theorem~\ref{surg-formula} follows from this.
In the rest of this section we give the details of the proof. 

\subsection{Variational formula} 
First we give a couple of lemmas which will help to establish the
result.

\begin{lemma}\label{traceclass} 
For $\mu > 0$, the operator $(\Lapsum + \mu)^{-1} - (\Lapplane + \mu)^{-1}$ is
trace class, the derivative of
$$
\log \det (\Lapsum + \mu)
$$
with respect to $\mu$ exists, and
\begin{equation}
\frac{d}{d\mu} \log \det (\Lapsum + \mu) = \tr \big( (\Lapsum + \mu)^{-1} - 
(\Lapplane + \mu)^{-1} \big).
\label{logdet}
\end{equation}
\end{lemma}

\begin{proof} See appendix.
\end{proof}

Next we need to introduce some notation. We define the Dirichlet and transmission
Poisson operators, $\Pdir(\mu)$ and $\Ptr(\mu)$, mapping from $H^{3/2}(H)$ to 
$H^2(\mathcal{O})
\oplus H^2(\Omega)$, respectively $H^{1/2}(H)$ to $H^2(\mathcal{O})
\oplus H^2(\Omega)$ by
\begin{equation}\begin{gathered}
\Pdir(\mu) (f) = u, \text{ where } (\Lapplus + \mu)u = 0 \text{ and } u \restriction
H = f \\
\Ptr(\mu) (f) = u, \text{ where } (\Lapplus + \mu)u = 0, \ u \text{ is continuous at }
H, \ [\pa_\nu u] = f.
\end{gathered}
\end{equation}
Here $[ ]$ denotes the jump in the argument at $H$ (the sign is
specified in Definition~\ref{Neumann-jump}. These are the Poisson operators
(more precisely, `half' of the Poisson operators) for the Dirichlet and transmission
boundary conditions considered by Forman. Notice that both map to the space
\begin{equation}
\{ u \in H^2(\mathcal{O}) \oplus H^2(\Omega) \mid u \text{ is continuous at } H \}.
\label{H2cts}
\end{equation}
Also note that both $(\Lapplane + \mu)^{-1}$ and $(\Lapplus + \mu)^{-1}$ map $L^2$
to \eqref{H2cts}. The corresponding trace operators, defined on \eqref{H2cts}, are
\begin{equation}\begin{gathered}
\Tdir (u) = u \restriction H \\
\Ttr (u) = [\pa_\nu u] .
\end{gathered}
\end{equation}
Thus, $R(\mu) = \Ttr \Pdir(\mu)$ and $R(\mu)^{-1} = \Tdir \Ptr(\mu)$. 

Then the following relations hold. 
\begin{lemma}
For $\mu > 0$, 
\begin{equation}
\frac{d}{d\mu} {\Pdir(\mu)} = -(\Lapplus + \mu)^{-1} \Pdir(\mu),
\label{Pdir-deriv}
\end{equation}
\begin{equation}
\Pdir(\mu) \Tdir \Ptr(\mu) = \Ptr(\mu),
\label{PTP}
\end{equation}
and
\begin{equation}
\Ptr(\mu) \Ttr (\Lapplus + \mu)^{-1} = (\Lapsum + \mu)^{-1} - (\Lapplane + \mu)^{-1}.
\label{second-form}
\end{equation}
\end{lemma}

\begin{proof} These are all routine. To prove \eqref{Pdir-deriv}, let $\Pdir(\mu) f
= u(\mu)$ and differentiate the equations
$$
(\Lapplus + \mu)u(\mu) = 0, \quad u(\mu) \restriction H = f
$$
and let $(d/d\mu) u = v$ to get
$$
(\Lapplus + \mu) v(\mu) = -u(\mu), \quad v(\mu) \restriction H = 0.
$$
It follows that 
$$
\big(\frac{d}{d\mu} \Pdir \big) f = v(\mu) = - (\Lapplus + \mu)^{-1} u(\mu) = 
- (\Lapplus + \mu)^{-1} \Pdir(\mu) f,
$$
which establishes \eqref{Pdir-deriv}. 

To prove the next formula, let $u = \Ptr(\mu) f $. Then $u$ is continuous at $H$; let
$g = \Tdir u$. Then $\Pdir(\mu) g$ is the unique solution $v$ of $(\Lapplus + \mu) v = 0$
which is continuous at $H$ and takes the value $g$ there. Thus $v=u$. This demonstates
\eqref{PTP}. 

To prove the final equation, consider the right hand side applied to $w$. This
gives us a function $u$ such that $(\Lapplus + \mu)u = 0$, with
$[\pa_\nu u]$ equal to $\Ttr (\Lapplus + \mu)^{-1} w$ 
(since $\Ttr (\Lapplane + \mu)^{-1} w = 0$). 
Therefore, $u = \Ptr (\Ttr (\Lapplus + \mu)^{-1} w)$, proving 
\eqref{second-form}.
\end{proof}

Finally we recall from \cite{BFK} that $(d/d\mu) R(\mu)$ is a
pseudodifferential operator of order $-1$ for $\mu > 0$, so
$R(\mu)^{-1} (d/d\mu) R(\mu)$ is an operator of order $-2$ and hence of
trace class. Thus, we can calculate 
\begin{alignat}{2}
-\frac{d}{d\mu} \big( \log \det R(\mu) \big) &= -\tr \big( R(\mu)^{-1} \frac{d}{d\mu}
{R(\mu)} \big) & \quad & \notag \\
&= -\tr \big( \Tdir \Ptr(\mu) \Ttr \frac{d}{d\mu}{ \Pdir(\mu) } \big) & \quad & 
\notag \\
&= \tr \big( \Tdir \Ptr(\mu) \Ttr (\Lapplus + \mu)^{-1} \Pdir(\mu) \big) & \quad 
&\text{
 by \eqref{Pdir-deriv} } \notag \\
&= \tr \big( \Pdir(\mu) \Tdir \Ptr(\mu) \Ttr (\Lapplus + \mu)^{-1} \big) & \quad & 
\notag \\
&= \tr \big( \Ptr(\mu) \Ttr (\Lapplus + \mu)^{-1} \big) & \quad
&\text{ by \eqref{PTP}} \notag \\
&= \tr \big( (\Lapplus + \mu)^{-1} - (\Lapplane + \mu)^{-1} \big)  & \quad &\text{ by 
\eqref{second-form}} \notag \\
&= \frac{d}{d\mu} \big( \log \det (\Lapplus + \mu) \big) & \quad
&\text{ by \eqref{logdet}}. \notag
\end{alignat}
This completes Step 1.

\vskip 20pt

\subsection{Asymptotics as $\mu$ tends to infinity} 
In this step we calculate the constant term in the asymptotic expansion
of the log determinants of $\Lapint + \mu$, $\Lapext + \mu$ and $R(\mu)$ as
$\mu \to \infty$. In fact, the interior Laplacian and the Neumann jump operator have
been treated in BFK, where it is shown that in both cases the constant term is zero, 
so we only need to deal with the exterior Laplacian. We use the
formula for the zeta function in terms of the regularized heat trace to deduce the
result. This method does not generalize very far, since it requires that the dependence on
$\mu$ is of the form $A + \mu$, but it has the advantage of being very explicit.

\begin{prop} The logarithm of the determinant of $\Lapext + \mu$ has an expansion
$$
\sum_{j = -2}^\infty \Big( p_j \mu^{- j/2} + q_j \mu^{- j/2} \log \mu \Big)
$$
as $\mu \to \infty$, with $p_0 = 0$. 
\end{prop}

\begin{proof} Recalling \eqref{zeta-mu-heat}, the zeta function for
$\Lapext + \mu$ is 
\begin{equation}
\zeta_{\Omega, \mu} (s) = 
\frac1{\Gamma(s)} \int_0^\infty t^s \rtr e^{-t\Lapext} e^{-\mu t} \frac{dt}{t}.
\label{zeta-mu-heattrace}
\end{equation}
This equation shows why there is an expansion as $\mu \to \infty$ with local
coefficients: the factor of $e^{-\mu t}$ means that the integral from $a$ to
infinity is exponentially decreasing in $\mu$, for any $a > 0$, 
so only the expansion of the
regularized heat trace at $t=0$
will contribute to polynomial-order asymptotics in $\mu$, and
this expansion is local. 

For precisely, for any integer $k$ we consider the expansion to $2k+2$ terms
of the regularized heat trace (see \eqref{ae}),
\begin{equation}
\sum_{j=-2}^{2k-1} a_j t^{j/2}
\label{finitesum}
\end{equation}
at $t=0$. 
Let $e_k(t)$ be the difference between the regularized heat trace of $e^{-t\Lapext}$ and
this finite expansion. Then, $e_k(t)$ is $O(t^k)$ as $t \to 0$ and, since the 
heat trace is bounded as $t \to \infty$, $e_k(t)$ is also $O(t^k)$ at infinity. It is
easy to see that if $e_k$ is substituted for the regularized heat trace in 
\eqref{zeta-mu-heattrace} then both the result, and the derivative in $s$ of the result,
is $O(\mu^{-k})$ as $\mu \to \infty$. Thus, to compute the expansion as 
$\mu \to \infty$ to this order we need only substitute \eqref{finitesum} in to 
\eqref{zeta-mu-heattrace}, namely
$$
\frac1{\Gamma(s)} \int_0^\infty t^s \sum_j a_j t^{j/2} e^{-\mu t} \frac{dt}{t}.
$$
Changing variable to $\overline{t} = t\mu$, this gives us
$$
\sum_j a_j \frac{\Gamma(s + j/2)}{\Gamma(s)} \mu^{-s-j/2} .
$$
Differentiating at $s=0$ gives us an expansion of the form above, with $p_0 = 0$
since the two $\Gamma$-factors cancel when $j=0$ to give a constant. 

\end{proof}

\subsection{Expansion as $\mu$ tends to zero}\label{mutozero} 
Here we consider the asymptotic expansion of the log determinants
\eqref{ext-asympt-0} --- \eqref{R-asympt-0}. The second of 
these, $\log \det (\Lapint  + \mu)$, is simply continuous as $\mu \to 0$, 
since $\Lapint + \mu$ has discrete
spectrum uniformly bounded away from zero as $\mu \to 0$. Thus 
\eqref{int-asympt-0} is obvious.

To understand the behaviour of $\log \det (\Lapext + \mu)$, recall from
\eqref{zeta-mu-chi} that the zeta function $\zeta_{\Omega, \mu}(s)$ is given by the
analytic continuation of 
\begin{equation}
\zeta_{\Omega, \mu} (s) = -\frac{s}{\pi} \int_0^\infty (\l^2 +
\mu)^{-s-1} \chi(\l)
\spl \l d\l + \frac1{\Gamma(s)} \int_0^\infty t^s e_2(t) e^{-\mu t} \frac{dt}{t}
\end{equation}
from $\Re s$ small, respectively $\Re s$ large. The contribution to
the log determinant from the second piece is continuous in $\mu$ as $\mu \to 0$,
so we get precisely $-\zeta_{\Omega, 2}'(0)$ (see \eqref{zeta-2-heat})
in the limit. 

In the first piece, the integrand is compactly supported, and
convergent uniformly near $s=0$ for fixed $\mu > 0$, since we have the
estimate \eqref{ltozero}. So the contribution to the log
determinant is equal to
\begin{equation}
\frac{1}{\pi} \int_0^\infty \frac{\l^2}{\l^2 + \mu} \spl 
\chi(\l) \frac{d\l}{\l}.
\label{chiterm}
\end{equation}
Let us write $\spl = \pi \ilg \l + \tilde{s}(\l)$, where, by Lemma~\ref{sc-phase-zero}, 
$\tilde{s}(\l)$ is
$O((\ilg \l)^2)$ as $\l \to 0$. (The function ilg is defined in
Definition~\ref{ilg}.) Replacing $s$ by $\tilde{s}$ in \eqref{chiterm}
makes the integral convergent uniformly down to $\mu = 0$, and we get a contribution
of 
\begin{equation}
\frac{1}{\pi} \int_0^\infty  \tilde{s}(\l) \chi(\l) \frac{d\l}{\l}.
\label{stilde}
\end{equation}

It remains to
consider what happens when $s$ is replaced by $\pi \ilg \l$. 
Thus, we are interested in
the expansion of the integral (where for convenience we replace
$\mu$ by $\nu^2$)
\begin{equation}
\int_0^\infty \frac{\l^2}{\l^2 + \nu^2} \ilg \l \, \chi(\l) \frac{d\l}{\l}.
\label{tough-int}
\end{equation}

To calculate this, we break up the integral into pieces. First consider
the integral from $0$ to $\nu$. We can estimate the absolute value by
$$
\int_0^\nu \frac{\l}{\nu^2} \ilg \l \, \chi(\l) \, d\l
= \nu^{-2} \int_0^\nu \l \ilg \l \, d\l.
$$
Since $(\l^2 \ilg \l/2)' \geq \l \ilg \l$ on the interval $[0, \nu]$, for
small $\nu$, this is estimated by
$$
\nu^{-2} \ [\frac{\l^2}{2} \ilg \l ]_0^\nu = O(\ilg \nu).
$$
Thus this term is $o(1)$ as $\nu \to 0$, and can be ignored.

Next consider the integral from $\nu$ to infinity of \eqref{tough-int}.
We claim that, up to an $O(\ilg \nu)$ error, we can replace the factor
$\l^2(\l^2 + \nu^2)^{-1}$ by $1$. To see this, we estimate the difference
$$
\int_\nu ^\infty \frac{\nu^2}{\l^2 + \nu^2} \ilg \l \, \chi(\l)\frac{d\l}{\l} \leq
\nu^2 \int_\nu ^2 \l^{-3} \ilg \l \, d\l. 
$$
Observe that, for small $\delta$, $\l^{-3} \ilg \l \leq \l^{-3} \ilg \l 
(2 - \ilg \l)$ on $[0, \delta]$, 
and the quantity on the right hand side is equal to the
derivative of $-\l^{-2} \ilg \l$. Therefore, this term is estimated by
$$
\nu^2 \ [(-\l^{-2} \ilg \l)]_\nu ^2 = O(\ilg \nu).
$$
Hence, up to $o(1)$ errors we are left with 
$$
\int_\nu ^\delta \ilg \l \frac{d\l}{\l} + \int_\delta^\infty \ilg \l
\chi(\l) \frac{d\l}{\l}.
$$
Let $\a = \ilg \l$, $\ep = \ilg \delta$ and $\td \chi(\a) = \chi(\l)$.
In these variables we have
\begin{equation}\begin{gathered}
\int_{\ilg \nu}^\ep \frac{d\a}{\a} + \int_\ep^\infty \td \chi(\a)
\frac{d\a}{\a} \\
= -\log \ilg \nu + \big( \int_\ep^\infty \td \chi(\a) \frac{d\a}{\a}
- \log \frac1{\ep} \big) \\
= \log \log \mu^{-1/2} + \HR \int_0^\infty \td \chi(\a) \frac{d\a}{\a},
\end{gathered}\end{equation}
where the last integral is a Hadamard regularized integral (see appendix).
This may be combined with \eqref{stilde} to give 
\begin{equation}
\log \det (\Lapext + \mu) =  \log \log \mu^{-1/2} +
\HR \int_0^\infty \td \chi(\a) \frac{d\a}{\a} - \zeta'_{\Omega, 2}(0)
+ \eqref{stilde} + o(1). 
\label{log-ext-form}
\end{equation}

\vskip 10pt

We need to compare this to the log determinant of $\Lapext$. By definition, this
is the coefficient of $-s$ in the expansion of the zeta function as $s \downto 0$. 
Recall that the zeta function is equal to $\zeta_{\Omega,1}(s) + \zeta_{\Omega,2}(s)$
as in \eqref{zeta-1-sp} and \eqref{zeta-2-heat}. Note that the contribution from
$\zeta_{\Omega, 2}$ is just $-\zeta'_{\Omega, 2}(0)$, matching one of the terms in
the expansion for 
$\Lapext + \mu$. From $\zeta_{\Omega, 1}$ we get a contribution which is 
the constant term, as $s \upto 0$, of
$$
\frac{1}{\pi} \int_0^\infty \l^{-2s} s(\l) \chi(\l) \frac{d\l}{\l}.
$$
Writing $s(\l) = \pi \ilg \l + \tilde{s}(\l)$ as before, with $s$ replaced by
$\tilde{s}$ above, the integral is convergent uniformly down to $s=0$ and
we get a contribution of exactly \eqref{stilde}. Thus, it remains to find the
constant term in the expansion of
\begin{equation}
\int_0^\infty \l^{-2s} \ilg \l \chi(\l) \frac{d\l}{\l}.
\label{r-int}
\end{equation}
Let us write $r = -s$, so $r \geq 0$. Let $\a = \ilg \l$; then
$\l^{-2s} = e^{-2r/\a}$ and 
$d\a/\a = \ilg \l \, d\l/\l$. Let $\ep > 0$ be arbitary. Then the
integral \eqref{r-int} is the same as 
$$ 
\int_\ep ^\infty e^{-2r/\a} \td \chi(\a) \frac{d\a}{\a} 
+ \int_{r/\ep} ^\ep e^{-2r/\a} \td \chi(\a) \frac{d\a}{\a} 
+ \int_0^{r/\ep} e^{-2r/\a} \td \chi(\a) \frac{d\a}{\a} .
$$
For small $\ep$, the factor $\td \chi(\a)$ may be replaced by one in the
last two integrals. Writing $\b = r/\a$, we get
$$
\int_\ep ^\infty e^{-2r/\a} \td \chi(\a) \frac{d\a}{\a} 
+ \int_{r/\ep} ^\ep e^{-2\b} \frac{d\b}{\b} 
+ \int_\ep ^\infty e^{-2\b} \frac{d\b}{\b} .
$$
It is clear that there is a divergent term $-\log r$ as $r \to
0$. We seek the limit when this divergent term is subtracted from
\eqref{r-int}. This limit is equal to 
$$ 
\int_\ep ^\infty \td \chi(\a) \frac{d\a}{\a} 
+ 2\log \ep - \int_0 ^\ep (e^{-2\b} - 1) \frac{d\b}{\b} 
+ \int_\ep ^\infty e^{-2\b} \frac{d\b}{\b} 
$$
for every $\ep > 0$. Taking the limit as $\ep \to 0$, the third term
disappears and one factor of $\log \ep$ combines with each of the
integrals to give two Hadamard-regularized integrals. Therefore, we
have shown that \eqref{r-int} has
an expansion
\begin{equation}
-\log r + \HR \int_0^\infty e^{-2\b} \frac{d\b}{\b}
+ \HR \int_0^\infty \td \chi(\a) \frac{d\a}{\a} + o(1), \ r \to 0.
\label{r-exp}\end{equation}
The first regularized integral appearing here
is equal to $\gamma + \log 2$ where $\gamma$ is Euler's constant
(\cite{AS}, chapter 5). The second regularized integral 
is the same one that appeared earlier. Combining it with \eqref{stilde}, 
we get the formula
$$
\log \det {}' \Lapext = \gamma + \log 2 +
\HR \int_0^\infty \td \chi(\a) \frac{d\a}{\a}  - \zeta'_{\Omega, 2}(0)
+ \eqref{stilde}.
$$
Comparing this with \eqref{log-ext-form}, we obtain \eqref{ext-asympt-0}.

\

Next we show \eqref{R-asympt-0}. We follow
the method of BFK. Since exactly one eigenvalue, say $\l_0(\mu)$, approaches
zero as $\mu \to 0$, we have
$$
\log \det R(\mu) = \log \l_0(\mu) + \log \det {} ' R + o(1)
$$
as $\mu \to 0$. Thus, we need to find the expansion of $\l_0(\mu)$; we use
BFK's characterization that
$$
\l_0(\mu)^{-1} = \text{ operator norm of } R(\mu)^{-1}.
$$

The operator $R(\nu^2)^{-1}$ is given by
$$
\omega \mapsto \Tdir \Ptr(\nu^2) \omega = \Tdir (\Lapplane + \nu^2)^{-1} J \omega.
$$
Here $J$ is the map $\omega \mapsto \omega \delta_H$, where $\delta_H$ is the
delta function supported on $H$. Let us analyze the behaviour as $\nu \to 0$.
By direct computation, it is not hard to show that
\begin{equation}
\frac1{|\xi|^2 + \nu^2} - 2\pi \log \frac1{\nu} \delta \quad \text{ converges in } 
\mathcal{S}'(\RR^2) \text{ as } \nu \to 0.
\label{conv}
\end{equation}
Moreover, away from the origin the convergence is as a symbol of order $-2$.
Conjugating with the Fourier transform yields an operator $M(\nu)
= (\Lapplane + \nu^2)^{-1} - 1/2\pi \log (1/\nu) 1$ whose kernel has a
pseudodifferential singularity (that is, conormal) at the diagonal, but with
growth as $|z-z'| \to \infty$. (Here, $1$ denotes the
operator whose kernel is identically equal to one, not the identity operator.)
Let $\rho$ be a function of compact support on $\RR^2$, which is identically equal
to one on a ball of large radius $R$. Then, since \eqref{conv} converges as a
symbol of order $-2$ away from the origin as $\nu \to 0$, $ \rho M(\nu) \rho$ converges
as $\nu \to 0$ as a pseudodifferential operator of order $-2$, and therefore,
as a bounded map from $H^{-1}(\RR^2)$ to $H^1(\RR^2)$. 

Therefore, 
\begin{equation}
R(\nu^2)^{-1} = \Tdir \rho \, (\Lapplane + \nu^2)^{-1} \rho \, J 
= \frac1{2\pi} \log \frac1{\nu} 1 + \Tdir \rho \, M(\nu) \rho \, J,
\end{equation}
where $1$ now denotes the operator on $L^2(H)$ with kernel equal to one.
Because $J$ is a continuous map from $L^2(H)$ to $H^{-1}(\RR^2)$ and $\Tdir$ is
a continuous map from $H^{1}(\RR^2)$ to $L^2(H)$, 
the second term is a family of operators on $L^2(H)$ with uniformly bounded
norm as $\nu \to 0$. Denoting the length of $H$ by $L$, the operator $1$ is
$L$ times a rank one projection on $L^2(H)$, so we see that the operator
norm of $R(\nu^2)^{-1}$ is equal to
$$
\| R(\nu^2)^{-1} \| = \frac{L}{2\pi} \log \frac1{\nu} + O(1)
$$
as $\nu \to 0$. 
Taking the logarithm of this, we see that
$$
\log (\l_0(\mu))^{-1} = \log \log \mu^{-1/2} + \log \frac{L}{2\pi} + o(1), 
\quad \mu \to 0.
$$
Thus, $\log \det R(\mu)$ has the expansion 
\begin{equation}
\log \det R(\mu) = -\log \log \mu^{-1/2} - \log \frac{L}{2\pi} + 
\log \det {}' R(0) + o(1), 
\quad \mu \to 0,
\end{equation}
which is \eqref{R-asympt-0}. This completes the proof of Theorem~\ref{surg-formula}.

\section{Compactness of Isophasal sets}

In section~\ref{isophasal} we defined the $\CI$ topology on exterior
domains. In this section we will prove
\begin{thm} \label{compactness}
Each class of isophasal planar domains is sequentially
compact in the $\CI$-topology.
\end{thm}

\subsection{Compactification of the problem}
First we explain how the problem is equivalent to a problem about a
bounded, non-flat domain on $S^2$. Consider an exterior domain
$\Omega$
whose closure does not contain the origin. 
Let $g$ be a metric on the plane of the form
$$
g = \frac{g_0(z)}{f(|z|^2)}, \quad z \in \RR^2,
$$
where $g_0$ is the flat metric, and $f(t)$ is equal to one for $t \leq R$
and $f(t) = t^2$ for $t \geq R'$. Then, $g$ is the same as the
standard metric on the ball $B(R, 0)$, whilst its behaviour at
infinity means that $g$ extends to a smooth metric on $S^2$ regarded
as the one-point compactification of $\RR^2$. This
compactifies $\Omega$ to a domain $\Omega '$ in $S^2$. Let $G$ be the
composition of inversion in the unit disc of $\RR^2$, followed by the
identification above to $S^2$. $G$ is then a conformal map from 
$\RR^2$ to $S^2$ which maps 
the inversion $\Omega^I$ of $\Omega$ to $\Omega '$. Thus
the metric $g$ on $\Omega'$ pulls back to the metric $e^{2\phi_0(w)} dw
d\overline{w}$ on $\Omega^I$,
where $w$ is a complex variable acting as
a coordinate on $\RR^2$ in the usual way. Here $\phi_0 = \log |G'|$; note that
$\phi_0$ is {\it not} harmonic, because $\Omega'$ is not flat. 
Let $F$ be a conformal map from
the unit disc $D$ to $\Omega^I$. Then the metric on $\Omega'$ pulls back
to $e^{2\phi_0(F(z))} e^{2\phi} dz d\overline z$, where $\phi = \log |F'|$
is a harmonic function (see \cite{OPS2}, section 1).

Now suppose that $\Omega$ varies within an isophasal class. Since the
topology on exterior domains is specified in terms of their
inversions, the first
thing we need to do is position $\Omega$ well with respect to the
inversion map. Recall that the first two heat invariants tell us that
the perimeter and area of $\Omega$ are fixed. By Lemma~\ref{inradius}
of the appendix, then, there is a uniform upper bound on the diameter,
and a uniform lower bound on the inradius, over the isophasal
class. Therefore, there is some $r > 0$ such that every $\Omega$ in
the isophasal class can be moved by a Euclidean motion so that the
boundary lies in the annulus
$$
A_r = \{ z \in \RR^2 \mid r < |z| < \frac1{r} \}.
$$
Consequently, the boundary of the inversion also lies inside $A_r$. We
choose $R$, in the definition of the metric $g$ above, to be larger
than $1/r$ so that all our domains $\Omega$ can be placed
isometrically within the flat part of the metric $g$. 

Next we compare our surgery formula for the exterior log determinant,
\begin{equation}
\log \det {}' \Lapext + \log \det \Lapint + \log \det R = \gamma +
\log \frac{L}{\pi}
\label{surg-1}\end{equation}
with
BFK's formula for $H \subset S^2$:
\begin{equation}
\log \det \Delta_{\Omega'} + \log \det \Lapint + \log \det R = \log \det 
\Delta_{(S^2,g)} + \log \frac{L}{A} .
\label{surg-2}\end{equation}
Notice that the two $R$ operators are the same, since the Laplacian is
conformally equivariant and thus the harmonic extension of a function on $H$ to
the exterior is the same whether we use the flat metric or the metric $g$ on the
exterior of the obstacle. The quantities $\log \det \Delta_{(S^2, g)}$ 
and $A$ are constant since we have fixed
$g$ once and for all, and as $\Omega$ varies over an isophasal class, the
log determinant of $\Lapext$ is fixed. 
Thus, subtracting the two equations we get
$$
\log \det \Delta_{\Omega'} = \text{ constant}
$$
over any isophasal class. In addition, since the metric on $S^2$ is
flat in a neighbourhood of the complement of $\Omega'$, the domains
$\Omega'$ all have the same heat invariants. 

Thus, our situation is that we have a class of metrics 
$e^{2\phi_0(F(z))} e^{2\phi} dz d\overline z$ on the unit disc, 
which have fixed determinant and heat invariants. This is the same
information as in \cite{OPS2}, but in this case the metric has an
extra factor of $e^{2\phi_0}$, where $\phi_0$ is evaluated at the
variable point $F(z)$. In the next two subsections we adapt the argument
of OPS to show that the set of $\phi$'s are
compact in $\CI$, and in the final subsection we use this to prove 
sequential compactness. 

{\it Remark. } We can express the log determinant of $\Lapext$
differently by considering also BFK's surgery formula for $S^2$ with
metric $g$ with respect to the unit disc, $D$. If we write $D'$ for the
complement of the unit disc, then this takes the form 
\begin{equation}
\log \det \Delta_{D} + \log \det \Delta_{(D', g)} + \log \det R_D = \log \det 
\Delta_{(S^2,g)} + \log \frac{2\pi}{A} .
\label{surg-3}\end{equation}
Noting that $\log \det \Delta_D$ and $\log \det R_D$ are universal
constants, we get by adding \eqref{surg-1} and \eqref{surg-3} and
subtracting \eqref{surg-2} that 
$$
\log \det {}' \Lapext = \log \det \Delta_{\Omega'} - \log \det
\Delta_{(D', g)} + \text{ universal constants }.
$$
Thus, up to universal constants, we can
compute the log determinant of an exterior domain $\Omega$ as a
difference of log determinants on two bounded domains. One is for the domain 
$\Omega'$ obtained by putting a metric
on the plane that is Euclidean on some large ball 
and conformally compactifies the plane at
infinity, and the other is the exterior of
the unit disc with respect to the same metric. It seems likely that
one could use this formula, together with a limiting process where $g$
becomes Euclidean on larger and larger balls, to find an explicit
Polyakov-Alvarez type formula for the exterior determinant.

\subsection{The Sobolev $\half$ estimate}
Let us begin by recalling the way in which OPS proved $\CI$ compactness
for isospectral planar domains.
Proving compactness is equivalent to obtaining uniform bounds on all
Sobolev norms of the function $\phi$ on the boundary of the disc $D$,
which determines the isometry class of the
domain $\obst$ as the image of a conformal map $F$. Specifically,
$\phi$ extends to a harmonic function on the disc, and then $F$ is the
unique analytic function satisfying $\log |F'| = \phi$, and $F(0) = 0$,
$F'(0) \geq 0$. Once uniform bounds on
the first few Sobolev norms of $\phi$ have been
obtained, it is easy to use Melrose's formulae for the heat invariants to
obtain uniform bounds on the higher Sobolev norms inductively. 

As in OPS, we need to place an additional constraint on $\phi$, namely
that $\phi$ is `balanced'. Since $F$ is only determined up to a
M\"obius transformation of the circle, $F$, and therefore $\phi$, are
not uniquely determined. OPS made the definition that $\phi$ is
balanced
if it satisfies 
$$
\intcirc e^{\phi} e^{i\theta} d\theta = 0.
$$
They proved that there is always a balanced conformal map from any
domain, flat or not, to the disc \cite{OPS2}. This condition is
important since then $\phi$ satisfies an improved inequality, as 
discussed below.

The crucial estimate of OPS is using the log determinant to obtain an 
$H^{1/2}$ estimate on $\phi$. This goes as follows: using the constancy
of the log determinant on the isospectral class, and the
Polyakov-Alvarez formula for the log determinant in terms of $\phi$, we obtain
\begin{equation}
\frac1{2} \int_{S^1} \phi \pa_n \phi + \int_{S^1} \phi = C.
\end{equation}
Here and below, $C$, $C_1$, etc, will denote constants which are uniform
over the isospectral class. The first term is almost equal to the square of
the Sobolev $1/2$-norm. In fact, expanding $\phi$ in a Fourier series,
$$
\phi(\theta) = \sum_n a_n e^{in\theta}, \quad a_{-n} = \overline{a_n},
$$
it is easy to show that 
$$
(2\pi)^{-1} \intcirc \phi \pa_n \phi \ d\theta = \sum_n |n| |a_n|^2,
$$
so
$$
\intcirc \phi \pa_n \phi \ d\theta + |a_0^2| = \intcirc \phi \pa_n \phi \ 
d\theta + \Big| \intcirc \phi \, d\theta \Big|^2
$$
is equivalent to the square of the Sobolev $\half$-norm. For future 
reference we note that 
\begin{equation}
\big\| \pa_\theta^j \pa_n \phi \big\|_{L^2} + \Big| \intcirc \phi \Big|
\label{Sob}
\end{equation}
is equivalent to the Sobolev $(j+1)$-norm. 

On the other hand, the Lebedev-Milin inequality
for balanced $\phi$ (\cite{OPS1}, equation (5) of the introduction) gives 
\begin{equation}
\log L - \intcirc \phi \leq \frac1{4} \intcirc \phi \pa_n \phi .
\label{Leb}
\end{equation}
Combining the two we find a bound on $\intcirc \phi \pa_n \phi$ and
then on $| \intcirc \phi |$, yielding a $H^{1/2}$ bound. (Without the
balanced hypothesis \eqref{Leb} is only valid with coefficient $\frac1{2}$
in front of the final term, which does not yield any Sobolev bound.) 

In our situation, the metric $e^{2\phi + 2\phi_0 \circ F} dz d\zbar$ is not
flat, so we get additional terms in our expression for the
log det. Let us write $\phit$ for $\phi + \phi_0 \circ F$. As above,
we may assume that $\phi$ is balanced. Also, since the area and length
are isophasal invariants, we have
\begin{equation}
\intcirc e^{\phit} = L, \quad \intdisc e^{2\phit} = A'
\label{length-vol}
\end{equation}
where $A' = \Area(S^2, g) - A$ and $A$ is the common volume of the obstacles
in our isophasal set. 

\begin{lemma} There is a uniform bound on the Sobolev half-norm of
both $\phi$ and $\phit$, regarded as functions on $S^1$,
as $\Omega$ ranges over an isophasal class. 
\end{lemma}

\begin{proof}
Constancy of the log determinant of $\Omega'$ over the isophasal
class implies that (\cite{OPS1}, equations (1.15), (1.16))
\begin{equation}
\intdisc |\nabla \phit|^2 + 2 \intcirc \phit + 3 \intcirc \pa_n \phit = C
\end{equation}
Expanding this out, we get
$$
\intdisc |\nabla \phi + \nabla (\phi_0 \circ F)|^2 +
2 \intcirc \phi + 2 \intcirc \phi_0 \circ F + 3 \intcirc \pa_n (\phi_0 \circ
F)  = C.
$$
Since $|a+b|^2 \geq 3/4 |a|^2 - 3 |b|^2$, we have
$$
\frac{3}{4}\intdisc |\nabla \phi|^2 - 3 \intdisc | \nabla (\phi_0 \circ F)|^2 +
2 \intcirc \phi + 2 \intcirc \phi_0 \circ F + 3 \intcirc \pa_n (\phi_0 \circ
F)  \leq C.
$$
Integrating the first term by parts, using sup bounds on $\phi_0$ and
$\nabla \phi_0$, and using $|F'| = e^\phi$, we get
$$
\frac{3}{4}\intdisc \phi \pa_n \phi - 3 C_1 \intdisc  e^{2\phi} +
2 \intcirc \phi + 2 C_2 + 3 C_3 \intcirc e^\phi  \leq C.
$$
Since we have \eqref{length-vol}, and 
$$
\intdisc e^{2\phi} = \intdisc e^{2\phit} e^{-2\phi_0 \circ F},
$$
the left hand side is bounded by $e^{2\| \phi_0 \|_\infty} A'$. The term
$\intcirc e^\phi$ is bounded similarly. Thus we have
\begin{equation}
\frac{3}{4}\intdisc \phi \pa_n \phi + 2 \intcirc \phi \leq C.
\end{equation}
Adding twice \eqref{Leb} to this inequality gives us a bound on
$\int \phi \pa_n \phi$, and then \eqref{Leb} provides a bound on
$| \int \phi |$. This gives a bound on the Sobolev one-half norm of
$\phi$. 

Next we bound the $H^{1/2}$ norm of $\phit$. Since 
$\phit = \phi + \phi_0 \circ
F$, it is sufficient to bound the $H^1$ norm of $\phi_0 \circ F$.
The derivative of $\phi_0 \circ F$ is bounded by
$| \phi_0' \circ F | |F'|$. We have uniform sup bounds on
$\phi_0'$. By Trudinger's inequality, (\cite{OPS1}, (3.63)),
\begin{equation}
\intcirc e^\phi \leq C \big\| \phi \big\|_{H^{1/2}},
\label{Trud}\end{equation}
so 
$$
\big\| |F'| \big\|_2^2 \leq \big\| e^\phi \big\|_2^2 \leq C 
\big\| \phi \big\|_{H^{1/2}}^2 \leq C,
$$
which gives us the required bound. 
\end{proof}

\subsection{Uniform $L^\infty$ and Sobolev $1$-bounds on $\phi$} 
To obtain a uniform bound on $|\phi|$,
we use the heat invariant $a_1$, which shows that \cite{OPS2}
$$
\int_H k^2(s) ds = \intcirc e^{-\phit} (1 + \pa_n \phit)^2 d\theta = C.
$$
Using $a^2 \leq 2(1+a)^2 + 2$, we obtain from this
$$
\intcirc e^{-\phit} (\pa_n \phit)^2 d\theta \leq 2C + 2\intcirc e^{-\phit}.
$$
But we have a uniform bound on $\| \phit \|_{H^{1/2}}$ so the Trudinger
inequality 
gives a uniform bound on $\int e^{-\phit}$. Thus,
$$
\intcirc e^{-\phit} (\pa_n \phit)^2 d\theta \leq C_1.
$$
Writing this in terms of $\phi$ and $\phi_0$, and using $a^2/2 - b^2 \leq 
(a+b)^2$, we obtain
$$
\intcirc e^{-\phi} \big( (\pa_n \phi)^2 - 2(\pa_n (\phi_0 \circ F))^2
\big) d\theta  \leq C_1.
$$
Since $|F'| = e^{\phi}$, we get
$$
\intcirc e^{-\phi} (\pa_n \phi)^2 \leq C_1 + 2\intcirc e^{\phi}
((\pa_n \phi_0) \circ F)^2 \leq C_2,
$$
applying the Trudinger inequality again. The argument of \cite{OPS2}, 
equation (1.13) on) can then be applied verbatim to conclude that
$\sup \phi$ and $\| \phi\|_{H^1}$ are uniformly bounded. 

\subsection{Higher Sobolev bounds} Here we will prove that for any
$k$, $\phi$ is uniformly
bounded in the Sobolev $k$-norm. The proof is by induction
on $k$. First we recall the proof in the OPS
case, where we just have the harmonic function $\phi$. 
Thus suppose that $\| \phi \|_{H^{j}}$ is uniformly bounded,
and, therefore, $\| \phi \|_{C^{j-1}}$ is also, by the Sobolev inequality. 
We start at $j=1$ since the hypotheses have been proved for this value
above. 

We use Melrose's result that
$$
\int_H (\pa_s^j k)^2(s) ds 
$$
is uniformly bounded. This implies that
$$
\intcirc (\pa_\theta^j k)^2 d\theta 
$$
is uniformly bounded, since 
$$
\left( \frac{d}{ds} \right)^j = \left( e^{-\phi} \frac{d}{d\theta} \right)^j
= e^{-j\phi} \left( \frac{d}{d\theta} \right)^j + \sum_{l=0}^{j-1}
p_l(e^{-\phi}, \pa_\theta \phi, \dots ,\pa_\theta^{j-1} \phi) \left(
\frac{d}{d\theta} \right)^l,
$$
where $p_l$ are polynomials, 
and all the occurrences of $\phi$ above are $L^\infty$-bounded by the
inductive assumption. 

Consider 
$$
\pa_\theta^j \pa_n \phi = \pa_\theta^j ( e^{\phi} k + 1);
$$
in view of \eqref{Sob}, it is sufficient to bound the $L^2$ norm of this
quantity. Taking derivatives, we find that this is equal to
\begin{equation}
\big( \pa_\theta^j e^{\phi} \big) k + \sum_{l=0}^{j-1} c_l \big( \pa_\theta^l
e^\phi \big) \big( \pa_\theta^{j-l} k \big).
\label{mess}
\end{equation}
The term $\pa_\theta^m e^{\phi}$ is equal to 
$$
(\pa_\theta^m \phi) e^{\phi} + q(e^{\phi}, 
\pa_\theta \phi, \dots ,\pa_\theta^{m-1} \phi),
$$
with $q$ a polynomial, and is therefore uniformly
$L^2$-bounded for $m=j$, and uniformly $L^\infty$-bounded for $m < j$. 
Thus, the $L^2$ norm of \eqref{mess} is bounded by
$$
C \Big( \| \pa_\theta^j e^{\phi} \|_2 \| k \|_\infty + \sum_l \| 
\pa_\theta^l e^\phi \|_\infty \| \pa_\theta^{j-l} k \|_2 \Big),
$$
which is a uniform bound. 

\

In the exterior domain case, uniform Sobolev
bounds on $\phit$ follow in the same way. Now we obtain uniform
Sobolev bounds on $\phi$. Since
\begin{equation}
\pa_\theta^j \pa_n \phit = \pa_\theta^j \pa_n \phi + \pa_\theta^j \pa_n 
(\phi_0 \circ F),
\label{phit-est}
\end{equation}
we need to obtain a uniform $L^2$-bound on the second term. Recall that
this depends on $\phi$ through the function $F$, since $\phi = \log |F'|$. 
Notice that
$$
\pa_n F = e^{i\theta} e^{\phi + i\psi},
$$
where $\psi$ is the harmonic conjugate of $\psi$. Thus
$$
\pa_\theta^j \pa_n (\phi_0 \circ F) = \pa_\theta^j \big(
(\pa_n \phi_0) \circ F (e^{i\theta} e^{\phi + i\psi}) \big)
$$
which involves at most $j$ derivatives of $\phi$ and $\psi$. Since $\psi$
is the harmonic conjugate of $\phi$, normal, resp. tangential 
derivatives of $\psi$ are equal to tangential resp. normal derivatives of
$\phi$ up to factors of $i$, and can therefore be estimated by derivatives
of $\phi$. The only way that $j$ derivatives of $\phi$ or $\psi$ can occur
is in the term
$$
\big( (\pa_n \phi_0) \circ F \big)
\big( \pa_\theta^j (e^{i\theta} e^{\phi + i\psi}) \big)
$$
which can be estimated by $\| \pa_n \phi_0 \|_\infty \| \pa_\theta^j
e^\phi \|_2$. 
The other terms can be estimated by $\| \phi_0 \|_{C^j} 
\| \pa_\theta^{j-1} e^\phi \|_{C^{j-1}}$. This shows that the second term 
of \eqref{phit-est}, and therefore also the first, is uniformly bounded
in $L^2$. This completes the inductive step of the proof. Thus, $\phi$
is uniformly bounded in $C^k$, for all $k$.

\subsection{Proof of Theorem~\ref{compactness}} Consider any sequence of
$\Omega_k$ of isophasal exterior domains. We will show that there is
a subsequence converging to some domain $\Omega$ in the same isophasal
class. 

First we show that there is a convergent subsequence. As discussed
above, we may assume that each $\Omega_k$ is placed so that its
boundary lies inside some annulus $A_r$. Let $\Omega_k^I$ be the
inversion of $\Omega_k$. Then we have shown above that any set of balanced
$\phi_k$ corresponding to the $\Omega_k^I$ 
is precompact in the $\CI$ topology. Consequently, the set of
normalized conformal maps $F_k$ (given by $\phi_k = \log |F_k'|$,
$F(0) = 0$, $F'(0) > 0$) is precompact in the $\CI$ topology
\cite{OPS2}. Passing to a subsequence, we may assume that the
$F_k$ converge in $\CI$ to a conformal map $F$. 

The
domain $\Omega_k^I$ is then given by $\Omega_k^I = E_k( F_k(D))$ for some
Euclidean motion $E_k$. Note that all $\Omega_k^I$ lie inside the ball
of radius $R$. Thus the collection $\{ E_k \}$ of isometries lie
inside some compact set; in fact, since $E_k(0)$ is contained in
$B_R(0)$, we have
$$
E_k \subset \{ T_b \circ R_\theta \mid \theta \in [0, 2\pi], \, |b|
\leq R \} \quad \forall \, k
$$
where $T_b$ is translation by complex number $b$ and $R_\theta$ is
rotation about the origin by angle $\theta$. Therefore, some
subsequence of $E_k$, say $E_{j_k}$, converge to a limiting Euclidean
motion $E$. Along this subsequence, $\td F_k = E_{j_k} \circ F_{j_k}$
converges in $\CI$, so the sequence of domains $\Omega_{j_k}^I$
converges to $\Omega^I \equiv E (F(D))$ . Thus, by definition of the
topology on exterior domains,
the sequence $\Omega_{j_k}$ converges to $\Omega$, the inversion of
$\Omega^I$. 
Thus the class of isophasal domains is precompact.

To complete the proof we have to
show that the isophasal class is closed in the $C^\infty$ topology; that
is, if a sequence of obstacles $\Omega_k$ with the same scattering phase
converges in $C^\infty$, then the limiting obstacle $\Omega$
has the same scattering phase. 

The scattering phase is given by
$$
s(\l) = -i \log \det S(\l) = -i \log \det \big( \Id + 
\sqrt{\frac{\l}{2\pi}} A(\l) \big),
$$
where $A(\l)$ is the operator on $L^2(S^1)$ with $C^\infty$ kernel
\begin{equation}
A(\l)(\theta, \omega) = \int_{|z| = R} u_\theta \frac{\pa e^{-i\l z 
\cdot \omega}}{\pa \nu}  - \frac{\pa u_\theta }{\pa \nu} e^{-i\l z \cdot \omega}
\label{amplitude}\end{equation}
for sufficiently large $R$ (see the appendix of 
\cite{HR}). Here $u_\theta$ is the distorted plane
wave with incoming direction $\theta$. Let us fix $\l$. 
The kernel $A(\l)$ 
belongs to the trace class of operators and $s(\l)$ is a continuous function
of $A(\l)$ with respect to the trace norm. Thus, we need to show that 
$A(\l)$ is a continuous function, as a trace class operator, of the
domain. It suffices to show that the distorted plane waves $u_\theta$,
together with a certain number of derivatives in $\theta$, are
continuous functions of the domains, say in $L^2_{\loc}$.
We will just treat the case of $u_\theta$ itself, since the
$\theta$-derivatives are treated similarly.

Choose a $\CI$ mapping $I_k : \overline{\Omega} \to \overline{\Omega_k}$
which is the identity for $|z| \geq 2R$. This can be done so that the
maps $I_k$ converge in $\CI$ to the identity map. Pulling back the
operator $\Lap - \l^2$ to $\Omega$ gives us differential operators
$L_k - \l^2$ on $\Omega$, with $L_k = \Lap$ for $|z| \geq 2R$, and
with coefficients converging to those of the Laplacian in $\CI$ as $k
\to \infty$. Thus, moving to $\Omega$, we must find solutions to
$$
(L_k - \l^2) u_k = 0, \quad u_k = 0 \text{ on } \pa \Omega,  \quad 
u_k - e^{-i\l z \cdot \theta} \text{ outgoing} 
$$
and show that $u_k \to u$, where $u$ solves the limiting problem on
$\Omega$. By standard methods this reduces to solving
$$
(L_k - \l^2) \td u_k = g_k, \quad \td u_k = 0 \text{ on } \pa \Omega,  \quad
\td u_k \text{ outgoing} 
$$
with $g_k \to g$ in $\CI$, and proving that $\td u_k \to \td u$, where
$\td u$
solves the limiting problem. These equations may be solved using the
method of the appendix of \cite{HR} (which in turn is based on the
proof of Theorem 5.2 of \cite{LP}). Tracing through the proof, it is
not hard to show that indeed $\td u_k \to \td u$ in $L^2_{\loc}$. 

Therefore, the scattering phase $s(\l)$ for $\Omega_k$ converges
pointwise to the scattering phase for $\Omega$. This shows that
$\Omega$ is isophasal with the $\Omega_k$ and completes the proof of
the theorem.

\section{Appendix}
\subsection{Uniform bounds on isophasal classes of domains}
In the proof of Theorem~\ref{compactness} we need the following lemma; 
recall that
any class of isophasal domains has fixed area and perimeter.

\begin{lemma}\label{inradius} 
For any obstacle $\obst$ there is a lower bound on the inradius,
and an upper bound on the circumradius, depending only upon the area and
perimeter of $\obst$. 
\end{lemma}

\begin{proof} The bound on the circumradius
$R \leq L$, where $L$ is the perimeter, is trivial. Consider the
inradius. Let $\obst$ be a domain with area $A$ and perimeter $L$ 
and $r$ a number such that $r/2$ is larger than the
inradius but smaller than the circumradius. 
It is possible to choose of covering of $\obst$ by balls 
of radius $5r$ whose centres lie in $\overline{\obst}$, such that
the balls of radius $r$ with the same centres are disjoint. To construct
such a covering, start with a finite covering by balls 
of radius $3r$ whose centres lie in $\overline{\obst}$. Choose any two
of the balls. If the balls of radius $r$ about their centres intersect, 
then the ball of radius $5r$ about any one of the centres contains both
balls of radius $3r$, and therefore one of the balls can be discarded.
By discarding balls successively in this way, we end up with a covering
with the required property. Let $N$ be the number of balls in the covering.

Considering the area of each ball, we have an inequality
\begin{equation}
25\pi r^2 N \geq A. 
\label{area}
\end{equation}
Since $r/2$ is assumed larger than the inradius, but smaller than the
circumradius, the circles of radius $r/2$ with the same centres as the
balls in our covering must all intersect the boundary. Since the circles
are all distance at least $r$ apart, this gives an inequality
\begin{equation}
L \geq N r.
\label{length}
\end{equation}
Combining the two inequalities, we find that
$$
25 \pi r \geq \frac{A}{L}.
$$
Taking the infimum over r, we find that
$$
\text{ inradius}(\obst) \geq \frac{2A}{25\pi L}.
$$
\end{proof}

\subsection{Computation of the scattering phase for the unit disc} To compute
the scattering phase, we construct the distorted plane waves for small $\l$.
These are given by
$$
\uol(z) = e^{-i\l z \cdot \omega} + \uol '(z), \quad z \in \RR^2, \ \omega
\in S^1
$$
where $\uol '$ satisfies
\begin{equation}\begin{aligned}
(\Delta + \l^2)\uol ' &= 0, \\
\uol ' (z) \restriction S^1 &= - e^{-i\l z \cdot \omega} \restriction S^1, \\
\uol' (z) &= |z|^{-1/2} e^{i\l |z|} a(\l, \hat z, -\omega) + O(|z|^{-3/2}),
\ |z| \to \infty .
\end{aligned}\end{equation}
Then the scattering phase is equal to $-i \log \det S(\l)$, where the
scattering matrix is the unitary operator
$$
S(\l) = \Id + \sqrt{\frac{\l}{2\pi}} A(\l),
$$
and $A(\l)$ has kernel (obtained by applying stationary phase to
\eqref{amplitude})
$$
A(\l)(\theta_1, \theta_2) = e^{i\pi/4} a(\l, \theta_1,
\theta_2).
$$
To find $\uol '$, we decompose $e^{-i\l z \cdot \omega} \restriction S^1$
as a Fourier series; $\uol '$ then has a corresponding decomposition
into Bessel functions. 

Let $\theta = \theta_1 - \theta_2$, where $\theta_i$ is now regarded as
a circular variable in $[0, 2\pi)$. By circular symmetry, $A(\l)$ depends
only on $\theta$. In terms of $\theta$, on the unit circle we have
\begin{equation*}\begin{aligned}
e^{-i\l z \cdot \omega} &= e^{i\l \cos \theta} \\
&= \sum_{j=0}^\infty \frac{(i\l)^j (e^{i\theta} + e^{-i\theta})^j}{2^j j!}.
\end{aligned}\end{equation*}
Hence, on the unit circle,
\begin{equation}\begin{aligned}
e^{-i\l z \cdot \omega} &= \sum_{n = -\infty}^\infty a_n(\l) e^{in\theta},
\text{ where } \\
a_0(\l) &= 1 + \l^2 \alpha_0(\l), \ \alpha_0(\l) \text{ bounded, } |\l| \leq 1,
\\
|a_n(\l)| &\leq \frac{\l^{|n|}}{|n|!}, \ |\l| \leq 1.
\label{an-est}
\end{aligned}\end{equation}
Then
$$
\uol ' (z) = \sum_n a_n(\l) \frac{H_{|n|}(\l |z|)}{H_{|n|}(\l)} e^{in \theta},
\theta = \hat z - (-\omega),
$$
where $H_n$ is the Hankel function of the first kind of order $n$. Using
large $|z|$ asymptotics of the Hankel function (\cite{AS}, chapter 9), this gives
$$
a(\l, \theta) = \sqrt{\frac{2}{\pi \l}} \frac{\sum_n a_n(\l)}{H_{|n|}(\l)}
e^{-i \pi/4}e^{-in \pi/2} e^{in \theta}.
$$
Since $|H_n(\l)| \leq \l^n$ for $\l \leq 1$, and we have the estimate
\eqref{an-est} for $a_n(\l)$, this series converges. As $\l \to 0$,
$$
H_0(\l) = \frac{2i}{\pi} \log \l + O(1), 
$$
so we obtain for $a(\l, \theta)$ the asymptotics
$$
a(\l, \theta) = i\sqrt{\frac{\pi}{2\l}} e^{i\pi/4} \ilg \l + O(\l^{-1/2}(\ilg \l)^2),
$$
and for $A(\l)$,
$$
A(\l)(\theta) = \frac{i}{2} \ilg \l + O((\ilg \l)^2),
$$
where the error term is a smooth function of $\theta$. Thus,
\begin{equation*}
s(\l) = -i \log \det S(\l) = -i \tr \log \big( \Id + A(\l) \big) = 
\pi \ilg \l + O((\ilg \l)^2).
\end{equation*}

\subsection{Proof of Lemma~\ref{traceclass}} The operator $(\Lapplus + \mu)^{-1}
- (\Lapplane + \mu)^{-1}$ is given in terms of the heat kernel by
\begin{equation}
\int_0^\infty e^{-t\mu} \big( e^{-t \Lapplus} - e^{-t \Lapplane} \big) dt.
\label{res}
\end{equation}
Jensen and Kato proved the following estimate of the trace norm of the
difference of heat kernels:
$$
\| \big( e^{-t \Lapplus} - e^{-t \Lapplane} \big) \|_1 = O(t^{-1/2}),
\quad t \to 0.
$$
On the other hand, if we denote the two semigroups by $H(t)$ and $H_0(t)$,
then
\begin{multline}
\| H(2t) - H_0(2t) \|_1 = \| H(t) \big( H(t) - H_0(t) \big) + 
\big( H(t) - H_0(t) \big) H_0(t) \|_1 \\
\leq \| H(t) \|_{\op} \| H(t) - H_0(t) \|_1 + \| H(t) - H_0(t) \|_1 \| H_0(t)
\|_{\op} = 2 \| H(t) - H_0(t) \|_1.
\end{multline}
Iterating this inequality shows that the trace of $H(t) - H_0(t)$ is $O(t)$
as $t \to \infty$. Thus, the integral \eqref{res} is convergent in trace
norm, and so the result is trace class. 

Next we prove the second part of the lemma. Using the functional calculus, we have
$$
\int\limits_{t}^\infty ( e^{-\tau(\Lapplus + \mu)} -
e^{-\tau(\Lapplane + \mu)} ) d\tau
= (\Lapplus + \mu)^{-1} e^{-t (\Lapplus + \mu)} - (\Lapplane + \mu)^{-1} 
e^{-t (\Lapplane + \mu)}.
$$
Using the estimates on trace norms above, we see that the right hand side 
is trace class, the trace is differentiable as a function of $t$,
and the derivative is minus the integrand on the left hand side of the equation. 
Hence we can calculate, for $\Re s > 1$,
\begin{equation}\begin{aligned}
{}&\frac{d}{d\mu} (\zeta_{\Omega, \mu}(s) + \zeta_{\obst, \mu}(s) ) 
= \frac{d}{d\mu} \frac1{\Gamma(s)}
\int_0^\infty t^s \tr \big( e^{-t(\Lapplus + \mu)} - e^{-t(\Lapplane + \mu)} \big)
\frac{dt}{t} \\
&= -\frac1{\Gamma(s)}
\int_0^\infty t^s \tr \big( e^{-t(\Lapplus + \mu)} - e^{-t(\Lapplane + \mu)} \big)
\, dt \\
&= \frac1{\Gamma(s)}
\int_0^\infty t^s \frac{d}{dt} \tr \Big( (\Lapplus + \mu)^{-1} e^{-t (\Lapplus + \mu)} 
- (\Lapplane + \mu)^{-1} e^{-t (\Lapplane + \mu)} \Big) \, dt \\
&= -\frac{s}{\Gamma(s)}
\int_0^\infty t^s \tr \Big( (\Lapplus + \mu)^{-1} e^{-t (\Lapplus + \mu)} 
- (\Lapplane + \mu)^{-1} e^{-t (\Lapplane + \mu)} \Big) \frac{dt}{t}.
\end{aligned}\end{equation}
The boundary term in the integration-by-parts is zero since, for $\Re s > 1$, 
the integrand tends to
zero at both zero and infinity. Since $(\Lapplus + \mu)^{-1}
- (\Lapplane + \mu)^{-1}$ is trace class, it follows that the integral has a simple
pole at $s=0$. But the term at the front, $s/\Gamma(s)$, has a double zero at
$s=0$, so when we take the minus derivative at $s=0$ to find the derivative of the log
determinant, we obtain precisely the pole of the integral at $s=0$. This pole is
$\tr ((\Lapplus + \mu)^{-1} - (\Lapplane + \mu)^{-1})$, proving \eqref{logdet}.
\qed

\subsection{Hadamard-regularized integrals and Pushforward}\label{pf}
Let $h(x)$ be smooth and compactly supported on $[0, \infty)$. Then
$\int h \, dx/x$ is not convergent. We define the Hadamard-regularized integral
of $h$ by the limit
$$
\HR \int_0^\infty h(x) \frac{dx}{x} = \lim_{\ep \to 0} \Big( \int_\ep ^\infty
h(x) \frac{dx}{x} - h(0) \log \frac1{\ep} \Big).
$$
It is easy to check that the limit exists. It may also be described as the
constant term in the asymptotic expansion of $\int_\ep ^\infty  h \, dx/x$ as
$\ep \to 0$.

Hadamard-regularized integrals turn up naturally in the Pushforward
theorem for polyhomogeneous functions proved by Melrose \cite{Mel4}, \cite{HMM}.
In fact, the authors first derived the expansions in section 3.3 using a
special case of this theorem, so we will include a brief discussion. 

The pushforward is invariantly defined on densities
rather than functions, so we consider densities defined on $\RR^2_+$. It is most natural
to consider b-densities, that is, densities of the form
$$
g(x_1, x_2) \frac{dx_1}{x_1} \frac{dx_2}{x_2}.
$$
One reason for this is that such densities have an invariantly defined
restriction to the boundary faces $x_i = 0$, obtained by cancelling the
factor $dx_i/x_i$ (which is invariant under changes of boundary defining
function at $x_i = 0$). Thus, the restriction of $g(x_1, x_2) \, dx_1/x_1
\, dx_2/x_2$ to
$x_1 = 0$ is $g(0, x_2) dx_2/x_2$. We will consider $g$ which are smooth and
have compact support. Then we have 

\begin{prop}\label{pushforward} 
Consider the map $f : \RR^2_+ \to \Rplus$ given by
$$
(x_1, x_2) \mapsto x = x_1 x_2.
$$
Let $u = v(x_1, x_2) \, dx_1/x_1 \, dx_2/x_2$ be a smooth b-density 
with compact support. Then the pushforward of $u$ has an asymptotic expansion
\begin{equation}
f_* u \sim \sum_{j=0}^\infty \big( p_j + q_j \log \frac1{x} \big) x^j 
\frac{dx}{x},
\label{exp-exists}
\end{equation}
where $q_0 = v(0,0)$ and 
\begin{equation}
p_0 = \big( \HR \int_0^\infty v(x_1, 0) \frac{dx_1}{x_1} +
\HR \int_0^\infty v(0, x_2) \frac{dx_2}{x_2} \big) .
\end{equation}
\end{prop}

This theorem is a special case of the general result in \cite{Mel4},
and is proved explicitly in \cite{HMM}, section 2. The expansion
\eqref{r-exp} follows immediately from this theorem by regarding the
integral as a pushforward under the map $(\a, \b) \mapsto r = \a\b$
(we also have to multiply by the formal factor $dr/r$ to make the
integrand into a density). In fact, the expansion \eqref{tough-int}
can also be deduced from the pushforward theorem by using the
operations of logarithmic and total blowup discussed in \cite{HMM}.

\bibliographystyle{plain}
\bibliography{sm}
\end{document}